\documentclass[11pt,a4paper,psamsfonts]{amsart}
\usepackage{amssymb,amscd,amsmath,amsthm}

\chardef\bslash=`\\ 

\newtheorem{theorem}{Theorem}[subsection]
\newtheorem{corollary}[theorem]{Corollary}
\newtheorem{lemma}[theorem]{Lemma}
\newtheorem{proposition}[theorem]{Proposition}

\theoremstyle{remark}
\newtheorem{remark}[theorem]{Remark}

\theoremstyle{definition}
\newtheorem{definition}[theorem]{Definition}

\numberwithin{equation}{subsection}

\newcommand{\thmref}[1]{Theorem~\ref{#1}}

\newcommand{\proref}[1]{Proposition~\ref{#1}}
\newcommand{\lemref}[1]{Lemma~\ref{#1}}

\newcommand{\N}{\mathbb N}

\newcommand{\Zp}{\mathbb Z_p}
\newcommand{\Q}{\mathbb Q} 

\newcommand{\primes}{\mathcal  P}
\newcommand{\nat}{\N^{\times}}

\newcommand{\qat}{\Q^*_+}

\newcommand{\leqr}{\preceq_r}
\newcommand{\alphatil}{\tilde{\alpha}}
\newcommand{\pitil}{\tilde{\pi}}
\newcommand{\Vtil}{\tilde{V}}
\newcommand{\Ad}{\operatorname{Ad}}

\title[dilations and full corners]{From endomorphisms to automorphisms and back:
dilations and full corners}

\date{January 25, 1998; revised October 12, 1998}
\thanks{Research supported by the Australian Research Council. \hfill Typeset by \AmS-\LaTeX}
\author[M. Laca]{Marcelo Laca}
\address{Department of Mathematics  \\
      The University of Newcastle\\  NSW  2308\\ AUSTRALIA}
\email{marcelo@math.newcastle.edu.au}

\subjclass{}
\begin{document}
\begin{abstract}

When $S$ is a discrete subsemigroup of a discrete group $G$ such that
$G = S^{-1} S$, it is possible to extend circle-valued multipliers
{}from $S$ to $G$; to dilate (projective) isometric representations
of $S$ to (projective) unitary representations of $G$; and to dilate/extend
actions of $S$ by injective endomorphisms of a C*-algebra 
to actions of $G$ by automorphisms of a larger C*-algebra. 
These dilations are unique provided they satisfy a minimality condition.
The (twisted) semigroup crossed product corresponding to an action  of $S$
is isomorphic to a full corner in the (twisted) crossed product by the 
dilated action of $G$.
This shows that crossed products by semigroup actions are
Morita equivalent to crossed products by group actions, making
powerful tools available to study their ideal structure and representation theory.
 The dilation of the system giving the Bost--Connes Hecke
C*-algebra from number theory is constructed explicitly as an application:
it is the crossed product $C_0(\mathbb A_f) \rtimes \qat$,
corresponding to the multiplicative action of the positive
rationals on the additive group $\mathbb A_f $ of finite adeles. 
\end{abstract}
\maketitle
\section*{Introduction}
In recent years there has been renewed interest in crossed products by
semigroups of endomorphisms, viewed now as universal algebras 
in contrast to their original presentation 
as corners in crossed products by groups.  
This new approach, initiated by 
Stacey \cite{sta} following a strategy pioneered by Raeburn 
for crossed products by group actions \cite{rae}, is based on 
the explicit formulation of a semigroup crossed product  
as the universal C*-algebra of a covariance relation.
As such, it motivated the development of specific techniques and
 brought about new insights and applications, e.g. 
\cite{sta,alnr,murnew,sri,twi-units,quasilat,bc-alg, hecke5,diri}. 
Nevertheless, the implicit view of semigroup crossed products as corners 
continues to have a very important role:
it is often invoked to prove the existence of nontrivial universal objects and 
it allows one to import results from the well-developed theory of crossed products by groups.
When the endomorphisms are injective and the semigroup 
is abelian the two approaches are equivalent,
and the proof involves using a direct limit to transform the 
endomorphisms into automorphisms and the isometries into unitaries. This has been done 
when the abelian semigroup is $\mathbb N$ \cite{cun,sta}, when it is totally ordered \cite{sri},
and, in general, when it is cancellative \cite{murnew}. 
As crossed products by more general (nonabelian) semigroups 
are being considered from the universal property point of view, 
the need arises to determine whether a realization
as corners in crossed products by groups is true and useful in those cases too.
This is the main task undertaken in the present work. 

A step away from commutativity of the acting semigroup was taken in \cite{semico} where
 isometric representations and multipliers of {\em normal} cancellative
semigroups were extended using the same direct limits 
(the semigroup $S$ is normal if $xS = Sx$ for every $x \in S$, 
in which case the natural notions of
 right and left orders on $S$ coincide).
Here we will go further and consider discrete semigroups that
 can be embedded
in a discrete group and for which the right order is cofinal;
since cofinality is a key ingredient of a directed system, this class is, arguably, 
the most general one for which the usual direct limit construction would work
without a major modification. 
 
Based on the results presented below one may argue that the relevant object  
is the action of an ordered group, and that there are two ways of looking at it;
the first is as an automorphic action on a C*-algebra
{\em taken together with a 
distinguished subalgebra which is invariant under the action of the positive cone},
and the second is simply as the endomorphic action of this positive cone
on the invariant subalgebra. We show that these two points of view are equivalent:
to go from the former to the latter 
one just cuts down the automorphisms to endomorphisms of the invariant subalgebra and
restricts to the positive cone, and the process 
is reversed by way of a dilation-extension construction,
\thmref{dil-ext}, which constitutes our first main result. 
We also  explicitly state and prove two
additional features of this equivalence that, in our opinion,
 have not previously received enough attention.
The first one is that the minimal automorphic dilation is canonically {\em unique},
which for instance allows one to test a good candidate, as done in Subsection \ref{diladeles} below.
The second one is that the crossed product by the semigroup action is 
realized as a {\em full} corner in the crossed product by a group action, so the equivalence of the 
two approaches technically translates into the strong Morita equivalence of the crossed products.
This is done in \thmref{fulcor}, which is our second main result.   

A modicum of extra work shows that these results are also valid for twisted crossed
products and projective isometric representations with circle-valued multipliers.
This requires the easy generalization, to Ore semigroups,
of results known for semigroups that are abelian \cite{arv,din,che,kle}
or normal \cite{semico,murpro}, which is done in the preliminary subsections \ref{multipls}
and \ref{twi.cross.prod}. 
The arguments given are for projective isometric representations and 
twisted crossed products, but setting all multipliers to be identically $1$ will lighten the burden 
slightly for those interested in the dilation-extension itself and not in 
projective representations, twisted crossed products, and extensions of multipliers.

In the final section we give an application to the semigroup dynamical system from number theory
\cite{bc-alg} which has the Bost-Connes Hecke C*-algebra \cite{bos-con} as its crossed product. 
Starting with the $p$-adic version of the system \cite[Section 5.4]{diri}
we show how one is quite naturally led to consider the ring of finite adeles 
with the multiplicative action of the positive rationals. 
This establishes a  natural heuristic link between the Bost-Connes Hecke C*-algebra
and the space $\mathcal A/\mathbb Q^*$, which lies at the heart
of Connes's recent formulation of the Riemann Hypothesis as a trace formula \cite{con-cr,con-rzf}. 

\section{Preliminaries}

In this first section we gather the basic  definitions and results concerning the 
semigroups on which we will be interested. We also generalize
other results about isometries and crossed products
that are valid, with more or less the same proofs, in the present setting,
although they were originally stated for particular cases.

\subsection{Ore semigroups.}

\begin{definition}
An {\em Ore semigroup} $S$ is a cancellative semigroup such that
 $Ss \cap St \neq \emptyset $ for every pair $s, t \in S$.
Ore semigroups are also known as {\em right--reversible} semigroups.
(We leave the obvious symmetric consideration of left--reversibility to the reader.)
\end{definition}

\begin{theorem}[Ore, Dubreil]
 A semigroup $S$ can be embedded in a group $G$ with $S^{-1} S = G$
if and only if it is an Ore semigroup.
In this case, the group $G$ is determined up to canonical isomorphism
and every semigroup homomorphism  $\phi$ from $S$ into a group $\mathcal G$ 
extends uniquely to a group homomorphism $\varphi : G \to \mathcal G$.  
\end{theorem}
\begin{proof}
See e.g. theorems 1.23, 1.24 and 1.25 in  \cite{cli-pre} for the first part. 
We only need to prove the assertion about extending $\phi$.
Since  $G = S^{-1} S$, given $x,y \in S$ there exist $u,v \in S$ such that
$v^{-1} u = y x^{-1}$, and hence the element $ux = vy $ is in $S x \cap S y$, proving 
that $S$ is directed by the relation defined by $s \leqr t$ if $t\in Ss$.
An easy argument shows that $\varphi(x^{-1} y) = \phi(x)^{-1} \phi(y) $ defines a
group homomorphism from $G = S^{-1} S $ to $\mathcal G$ that extends $\phi$.
\end{proof}
\begin{remark}
The last assertion of the theorem generalizes \cite[Lemma 1.1]{semico}.
Here we have found it more convenient, for compatibility with the rest of \cite{semico},
to work with the {\em right order} $\leqr$  determined by $S$ on $G$ via
$x \leqr y$ if $y \in S x$.
\end{remark}

To illustrate the class  of semigroups being considered, we list a few examples
which have appeared recently in the context of semigroup actions: 
\begin{itemize}
\item Abelian semigroups, (notably the multiplicative nonzero integers in an algebraic number field
\cite{hecke5});

\item Semigroups obtained by pulling back the positive cone from a totally ordered quotient \cite{phi-rae};

\item Normal semigroups, in particular semidirect products \cite{semico,murpro}; 

\item  Groups of matrices over the integers having positive determinant 
\cite[Example 4.3]{bre};

\end{itemize}

\subsection{Extending multipliers and dilating isometries.}\label{multipls}

Let $\lambda$ be a circle--valued multiplier on $S$, that is, a function
$\lambda : S\times S \to \mathbb T$  such that 
$$
\lambda(r,s) \lambda(rs,t) = \lambda(r,st) \lambda(s,t), \quad r,s,t \in S.
$$
A {\em projective isometric representation} of $S$ with multiplier $\lambda$
on a Hilbert space $H$
(an isometric $\lambda$--representation of $S$ on $H$)
is a family $\{V_s: s\in S\}$ of isometries on $H$
such that $V_s V_t = \lambda(s,t) V_{st}$.

A twisted version of Ito's dilation theorem  \cite{ito} was obtained in 
\cite[Theorem 2.1]{semico}, where
projective isometric representations of normal semigroups were dilated
to projective unitary representations.
Essentially the same proof, inspired on Douglas's \cite{dou},
works for Ore semigroups and gives the following.

\begin{theorem}\label{dilation}
Suppose $S$ is an Ore semigroup 
 and let $\{V_s: s \in S\}$ be an isometric $\lambda$--representation
of $S$ on a Hilbert space $H$, where $\lambda$ is a multiplier on $S$.
Then there exists a unitary $\lambda$--representation of $S$ on a Hilbert space $\mathcal H$
containing a copy of $H$ such that
\begin{enumerate}
\item[(i)] $U_s$ leaves $H$ invariant and $U_s |_H = V_s$; and
\item[(ii)] $\bigcup_{s\in S} U_s^*H$ is dense in $\mathcal H$.
\end{enumerate}
\end{theorem}
\begin{proof}
 Verbatim from the proof of \cite[Theorem 2.1]{semico},
except for the following minor modification of the part of the argument 
where normality is used to obtain an
admissible value for the function $f_t $. 
The value $st$ used there has to be substituted by any (fixed) $z \in Ss \cap St$,
and thus the fourth paragraph there should be replaced by the following one.

Suppose now that $f \in H_0$ and $t \in S$, and consider 
the function $f_t $ defined by 
$f_t(x) = \lambda(x,t) f(xt)$ for $x \in S$. If $s \in S$ is admissible for $f$,
let $z \in  Ss \cap S t$. We will show that
$s_0 := zt^{-1}$ is admissible for $f_t$.  For every $x \in Ss_0$, $xt \in  Sz$, 
and since $z$ is admissible for $f$
\begin{eqnarray*}
\lambda(x,t)f(xt) &= &\lambda(x,t) \overline{\lambda(xtz^{-1}, z)}V_{xtz^{-1}} f(z)\\
        & = & \overline{\lambda(xtz^{-1}, zt^{-1})}V_{xtz^{-1}}\lambda(zt^{-1},t)f(zt^{-1}t)\\
        & = &  \overline{\lambda(xs_0^{-1}, s_0)}V_{xs_0^{-1}} f_t(s_0)
\end{eqnarray*}
where the second equality holds by the multiplier property applied to the elements 
$xtz^{-1}$, $zt^{-1}$, and $t$ in $S$.
This proves that $s_0$ is admissible for $f_t$, so $f_t \in H_0$.
\end{proof}

Since the  results of \cite{semico} concerning discrete normal semigroups
depend only on this dilation theorem and on the unique
extension of group--valued homomorphisms,
they too are valid for Ore semigroups and we list them here
for reference. 

\begin{theorem}\label{semimult}
Suppose $S$ is an Ore semigroup and let $G = S^{-1} S$.  Then 
\begin{enumerate}
\item Every multiplier on $S$ extends to a multiplier on $G$.
\item Restriction of multipliers on $G$ to multipliers on $S$ gives an
isomorphism of $H^2(G,\mathbb T)$ onto  $H^2(S,\mathbb T)$.
\item Suppose $\lambda$ is a multiplier on $S$ and let $V$ be
 a $\lambda$--representation of $S$ by isometries on $H$.
Assume $\mu$ is a multiplier on $G$ extending $\lambda$. 
Then there exists a
unitary $\mu$--representation $U$ of $G$ on a Hilbert space $\mathcal H$
 containing a copy of $H$ such that
$U_s|H = V_s$ for $s \in S$, and $\bigcup_{s\in S}U_s^* H$ dense in $\mathcal H$.
Moreover, $U$ and $\mathcal H$ are unique up to canonical isomorphism.
\end{enumerate}
\end{theorem} 
\begin{proof}
The proofs of all but the last statement about uniqueness are as in 
Theorem 2.2, Corollary 2.3 and Corollary 2.4 of \cite{semico}, provided
one considers the left-quotients $x = t^{-1} s$ 
instead of the right-quotients used there. 
 In order to prove the uniqueness statement suppose
$(U', \mathcal H ')$ is another  unitary $\mu$-representation
such that ${U'}_s|H = V_s$ and $\bigcup_{s\in S}{U'}_s^* H$ is dense in $\mathcal H'$.
It is easy to see that the map 
$$
W: U_s^* h \mapsto {U'}_s^* h , \qquad  s\in S,  h \in H
$$ 
is isometric, and that it extends to an isomorphism
of $\mathcal H $ to $ \mathcal H'$ because of the density condition. 
It only remains to show that $W$ intertwines $U$ and $U'$.
Since $S$ is an Ore semigroup, for every $x$ and $s$ in $S$ there exist $z $ and $t$ in $S$ 
such that $x s^{-1} = t ^{-1} z$.  Then $tx = zs$, so
\begin{eqnarray*}
W U_x (U_s^* h) &=& W U_x U_{tx}^* U_{zs} U_s^* h = 
 \mu(t,x)  \overline{\mu(z,s)} W U_t^* U_{z} h =  \mu(t,x)  \overline{\mu(z,s)} W U_t^* (V_z h)\\
                &=& \mu(t,x)  \overline{\mu(z,s)} {U'}_t^* (V_z h) =  \mu(t,x)  \overline{\mu(z,s)}{U'}_t^* {U'}_z  h  
= {U'}_x  ({U'}_s^*  h ) \\
&=& {U'}_x W  (U_s^* h)
\end{eqnarray*}
This shows that $WU_x = U'_x W$ for every $x\in S$,  hence for every $x \in G$.
\end{proof}

\subsection{Twisted semigroup crossed products}\label{twi.cross.prod}
 Suppose $A$ is a unital C*-algebra and let $\alpha$ be an action of the discrete semigroup $S$ by 
not necessarily unital endomorphisms of $A$.
Let $\lambda$ be a circle-valued multiplier on $S$.
A {\em twisted covariant representation} 
of the  semigroup dynamical system $(A, S, \alpha)$
with multiplier $\lambda$ is a pair $(\pi, V)$ in which
\begin{enumerate}
\item $\pi$ is a unital representation of $A$ on $H$,
\item  $V: S \to Isom( H)$  is a projective isometric representation of
 $S$ with multiplier $\lambda$, i.e., $V_s V_t = \lambda(s,t) V_{st}$,
and
\item the covariance condition 
$\pi(\alpha_t(a)) = V_t \pi(a) V_t^*$ 
holds for every $a\in A$ and $ t\in S$.
\end{enumerate}
When dealing with twisted covariant
representations with a specific multiplier $\lambda$, we will refer to the dynamical
system  as a twisted dynamical system and
denote it by $(A, S, \alpha, \lambda)$.

The (twisted) crossed product associated to $(A, S, \alpha, \lambda)$ is a
 C*-algebra $A\rtimes_{\alpha, \lambda}S$  
together with a unital homomorphism 
$i_A :A \rightarrow A\rtimes_{\alpha, \lambda}S$
 and a projective $\lambda$-representation of $S$ as isometries 
$i_{S}: S \rightarrow A\rtimes_{\alpha, \lambda}S
$ 
such that
\begin{enumerate}
\item $(i_A, i_{S})$ is a twisted covariant representation for
$(A, S, \alpha, \lambda)$,
\item for any other covariant representation $(\pi, V)$ there is a representation
$\pi \times V$ of $A\rtimes_{\alpha, \lambda}S$ such that
$\pi = (\pi\times V )\circ i_A$ and $V = (\pi\times V)
 \circ i_{S}$, and
\item  $A\rtimes_{\alpha, \lambda}S$ is generated by $i_A(A)$ and
$i_{S} (S)$ as a C*-algebra.
\end{enumerate}
The existence of a nontrivial universal object 
associated to $(A, S, \alpha, \lambda)$ depends 
on the existence of a nontrivial twisted covariant representation
with multiplier $\lambda$. For general endomorphisms
such representations need not exist, even in the untwisted case.
For instance, the action of $\mathbb N$ by surjective shift-endomorphisms of $c_0$ described
in Example 2.1(a) of \cite{sta} does not admit any nontrivial covariant representations.
We will assume that our endomorphisms are injective, hence
nontriviality of the semigroup crossed product will follow from its
realization as a corner in a nontrivial classical crossed product. 
See \cite{sta,murnew} for abelian semigroups, and Remark \ref{nontriv} below. 
There are other possible covariance conditions which
yield nontrivial crossed products 
even if the endomorphisms fail to be injective,
see e.g. \cite{murpac} and \cite{pet}.
We will not deal with them here, but we refer 
to \cite{lam} for an interesting comparative discussion of  
the different constructions.

\begin{remark}
It is immediate from the definition that
the crossed product $A \rtimes S$ is generated, as a C*-algebra, by the monomials
$v_x^* a v_y $ with  $a \in A$, and $ x,y \in S$, but more is true for Ore semigroups:
 the products of such monomials 
can be simplified using covariance to obtain another monomial of the same type.
Specifically, in order to simplify the product $v_x^* a v_y v_r^* b v_s$ we begin by
finding elements $t$ and $z$ in $S$ such that $y r^{-1}= t^{-1} z$, so that
$ty = zr$, (such elements do exist because $S$ is an Ore semigroup). 
It follows that
\begin{eqnarray*} 
  v_x^* a v_y v_r^* b v_s 
&=& \lambda(y,t) \overline{\lambda(z,r)} v_x^* a v_y v_{ty}^* v_{zr} v_r^* b v_s\\
&=& \lambda(y,t) \overline{\lambda(z,r)} v_x^* a v_yv_y^* v_t^*  v_z v_rv_r^* b v_s\\
&=& \lambda(y,t) \overline{\lambda(z,r)} v_x^* v_t^* 
    \alpha_t(a\alpha_y(1)) \alpha_z(\alpha_r(1) b ) v_z v_s\\
&=& \lambda(y,t) \overline{\lambda(z,r)}  \overline{\lambda(t,x)} \lambda(z,s)
     v_{tx}^*  \alpha_t(a\alpha_y(1)) \alpha_z(\alpha_r(1) b ) v_{zs},
\end{eqnarray*}  
hence the linear span of such monomials is dense in the crossed product.
\end{remark}

\section{The minimal automorphic dilation.}
\label{min-aut-ext}
There are two steps in realizing
a semigroup crossed product as a corner in a crossed
product by a group action. The first one
is the dilation-extension of a semigroup action by injective endomorphisms to a group action
by automorphisms,
and the second one is the corresponding dilation-extension of covariant representations
of the semigroup dynamical system to covariant representations
of the dilated system.

\subsection{A dilation-extension theorem.}
\begin{theorem}\label{dil-ext}
Assume $S$ is an Ore semigroup with enveloping group $G = S^{-1} S$ and
let $\alpha$ be an action of $S$ by injective endomorphisms of a unital C*-algebra $A$. 
Then there exists a C*-dynamical system $(B, G, \beta)$, 
unique up to isomorphism, consisting of
an action $\beta$ of $G$ by automorphisms of a C*-algebra $B$
and an embedding $i: A \to B$ such that
\begin{enumerate}
\item $\beta$ dilates $\alpha$,  that is, 
$\beta_s\circ i = i \circ \alpha_s$ for $s \in S$, and
\item $(B, G, \beta)$ is minimal,  that is, $\bigcup_{s\in S}\beta_s^{-1}(i(A))$ is dense in $B$.
\end{enumerate} 
\end{theorem}

\begin{proof}
By right--reversibility, $S$ is directed by $\leqr$ so
one may follow the argument of \cite[Section 2]{murnew}. 
However, extra work is needed here: 
since $G$ need not be abelian, the choice of embeddings 
in the directed system must be carefully matched to the choice of
right-order $\leqr$ on $S$.

Consider the directed system of C*-algebras determined 
by the maps 
$\alpha_y^x = \alpha_{yx^{-1}}$ from $A_x := A $ into $A_y := A$, 
for $x \in S$ and $y \in  Sx$,
 i.e. for $x \leqr y$ in $S$.
By \cite[Proposition 11.4.1(i)]{kad-rin} there exists an inductive
limit C*-algebra
$A_\infty$ together with embeddings $\alpha^x : A_x \to A_\infty$
 such that $\alpha^x = \alpha^y \circ \alpha^x_y$ whenever $x\leqr y$,
and such that $\bigcup_{x\in S} \alpha^x(A_x)$ is dense in $A_\infty$.

The next step is to extend the endomorphism  $\alpha_s$ to an automorphism of $A_\infty$.
For any fixed $s \in S$ the subset $Ss$ of $S$ is cofinal, so $A_\infty$ is also the 
inductive limit of the directed subsystem $(A_x, x\in Ss)$,
and, for this subsystem, we may consider new embeddings
$\psi^x : A_x \to A_\infty$ defined by $\psi^x (a) = \alpha^{xs^{-1}}(a)$ 
for $x \in Ss$ and $a \in A_x$.
By \cite[Proposition 11.4.1(ii)]{kad-rin} there is an automorphism
$\alphatil_s$ of $A_\infty$ such that $\alphatil_s \circ \alpha^x = \psi^x$  for every $x \in Ss$.
Since $  \alpha^1 =  \alpha^s \circ \alpha^1_s$ and $\psi^x = \alpha^{xs^{-1}}$,
the choice $x = s$ gives
$$
\alphatil_s \circ \alpha^1 
=  \alphatil_s \circ \alpha^s \circ \alpha^1_s = \alpha^1 \circ \alpha_s
$$
so that (1) holds with $\beta = \alphatil$ and $i = \alpha^1 : A_1 \to A_\infty$.

Since $\alphatil_s^{-1}(i(A)) = \alpha^s (A_s)$, (2) also holds. Uniqueness
of the dilated system follows from \cite[Proposition 11.4.1(ii)]{kad-rin}:
 $A_\infty$ is the closure of the union of the subalgebras 
$\alphatil_s^{-1}(i(A))$ with $s \in S$, if 
$(B, G, \beta)$ is another minimal dilation with embedding $j : A \to B$
then there is an isomorphism 
$\theta: A_\infty \rightarrow  B$ given by 
$\theta \circ \alphatil_{s^{-1}}( i(a)) = \beta_{s^{-1}} (j(a))$
 for $a \in A$ and hence which intertwines $\alphatil$ and $\beta$. 
\end{proof}

\begin{definition}
A system $(B, G, \beta)$ satisfying the conditions (1) and (2) of \thmref{dil-ext}
is called the {\em minimal automorphic dilation} of  $(A,S,\alpha)$. 
If $\lambda$ is
a multiplier on $S$ with extension $\mu$ to $G$, we say that 
the twisted system $(B, G, \beta, \mu)$
is the  minimal automorphic dilation of the twisted system
$(A,S,\alpha,\lambda)$. (By \thmref{semimult} the extended multiplier exists
and is unique up to a coboundary.)
\end{definition}

\begin{lemma}\label{dil-cov-rep}
Let $(\pi,V)$ be a covariant representation for the twisted system 
$(A,S,\alpha,\lambda)$ on the Hilbert space $H$, 
and let $\Vtil$ be the minimal projective unitary dilation of $V$ on $\mathcal H$
given by \thmref{dilation}.
 Then there exists a representation $\pitil$ of $B$ on $\mathcal H$ such that 
$(\pitil,\Vtil)$ is covariant for the minimal automorphic dilation $(B, G, \beta, \mu)$
and $\pitil  \circ i = \pi$ on $H$.
\end{lemma}
\begin{proof}
We work with the dense subspace $\mathcal H_0 = \bigcup_{t\in S} U_t^* H$ of $\mathcal H$
and the dense subalgebra $B_0 = \bigcup_{s\in S} \beta_s^{-1}(i(A))$.
If $\xi \in \mathcal H_0$ there exists $t \in S$ such that $U_t \xi \in H$;
 assume $b = \beta_t^{-1}(i(a))$, since we want $(\pitil, \Vtil)$ to be covariant,
the only choice is to define $\pitil$ by
$$
\pitil(b) \xi = \pitil(\beta_t^{-1}(i(a))) \xi 
= \Vtil_t^* \pitil(i(a)) \Vtil_t \xi = \Vtil_t^* \pi(a) \Vtil_t \xi
$$ 
because $\pitil$ restricted to $i(A)$ and cut down to $H$ has to be equal to $\pi$.

Of course we have to show that this actually defines an operator 
$\pitil (b)$ on $\mathcal H$ for each $b \in B_0$, 
that $\pitil$ extends to a homomorphism from all of $B$ to $B(\mathcal H)$,
and that $(\pitil,\Vtil)$ is covariant.

The first step is to  define $\pitil(b)$ on $\mathcal H_0$ for a fixed $b \in B_0$.
We begin by fixing $b \in B_0$, $a \in A$ and $s\in S$ such that
 $b = \beta_s^{-1}(i(a))$. For $\xi \in \Vtil^*_{t_0} H$
 with $t_0$ in the cofinal set $Ss$, we let
\begin{equation} \label{dil-rep}
\varphi(b) \xi = \Vtil_{t_0}^* \pi (\alpha_{{t_0}s^{-1}}(a)) \Vtil_{t_0} \xi.
\end{equation}
If $t \in St_0$ then $\xi \in \Vtil^*_t H$, 
 and
\begin{eqnarray*}
\Vtil_{t}^* \pi (\alpha_{{t}s^{-1}}(a)) \Vtil_{t} \xi_0  &=&
\Vtil_{t_0}^* \Vtil_{tt_0^{-1}}^* \pi(\alpha_{t t_0^{-1}} \circ \alpha_{{t_0}s^{-1}}(a)) 
\Vtil_{tt_0^{-1}} \Vtil_{t_0} \xi  \\
& = & \Vtil_{t_0}^*    \Vtil_{tt_0^{-1}}^* V_{t t_0^{-1}} 
                                     \pi( \alpha_{{t_0}s^{-1}}(a))
                       V_{t t_0^{-1}}^* \Vtil_{tt_0^{-1}}       
  \Vtil_{t_0} \xi \\
& = & \Vtil_{t_0}^*  \pi( \alpha_{{t_0}s^{-1}}(a))  \Vtil_{t_0} \xi. 
\end{eqnarray*}
So the definition of $\phi (b) \xi $ could have been given using any  $t \in St_0$ 
in place of $t_0$. Next we show that $\phi (b) \xi $ is also 
independent of $s$ and $a$, in the sense that
if $b$ is also equal to $ \beta_{s'}^{-1}(i(a')) $ then 
$\alpha_{{t}{s'}^{-1}}(a')$ is equal to $\alpha_{{t} s^{-1}}(a)$
for $t$ in a cofinal set.
To see this let $t \in Ss \cap Ss'$. 
Then $\alpha^t \circ \alpha^{s'}_t (a') = \alpha^{s'}(a') = \beta_{s'}^{-1}(i(a')) =
\beta_{s}^{-1}(i(a)) =  \alpha^{s}(a) = \alpha^t \circ \alpha^{s}_t (a)$,
and since the embedding $\alpha^t$ is injective, it follows that
$\alpha_{t{s'}^{-1}}(a') = \alpha_{ts^{-1}}(a)$.

The map $\varphi(b) : \mathcal H_0 \to \mathcal H_0$ is clearly linear, and 
since the endomorphisms are injective, $\| \varphi (b) \xi\| \leq \|b\| \| \xi \|$. 
Thus $\phi(b)$ can be uniquely extended to a bounded linear operator 
(also denoted $\varphi(b)$) 
on all of $\mathcal H$ such that $\| \varphi(b)\| \leq \|b\|$.
For any $s$ the map
$\Ad_{\Vtil_{t_0}^*} \circ \pi \circ \alpha_{{t_0}s^{-1}}$ is a *-homomorphism 
on $A$, and   by cofinality of $\leqr$, for any  $b_1$ and $ b_2 $ in $B_0$  
there exist $s \in S$ and $a_1$ and $a_2$ in $A$ such that 
$b_1 = \beta_s^{-1}(i(a_1))$ and $b_2 = \beta_s^{-1}(i(a_2))$.
It follows easily from (\ref{dil-rep}) that $\varphi : B_0 \to B(\mathcal H)$
is a *-homomorphism which can be extended to 
a representation $\pitil$ of $B$ on $\mathcal H$.

Putting $a = 1$ in (\ref{dil-rep}) shows that $\pitil$ is nondegenerate and there
only remains to check that
 $(\pitil,\Vtil)$ is a covariant pair for $(B,G, \beta, \mu)$.

Suppose first $x \in S$ and $b \in B_0$; we can assume that $b = \beta^{-1}_s (i(a))$
for some $a \in A$ and $s \in Sx$.
Let  $\xi \in \Vtil^*_t H$; we can assume $t \in Ss \subset Sx$, and
we observe that $\Vtil_x \xi \in \Vtil^*_{tx^{-1}} H$. Then
\begin{eqnarray*}
\pitil(\beta_x(b)) \Vtil_x \xi & = & \pitil(\beta_{x s^{-1}}(i(a))) \Vtil_x \xi\\
    & = &   \pitil(\beta^{-1}_{sx^{-1}}(i(a))) \Vtil_x \xi\\
    & = & \Vtil^*_{tx^{-1}} \pi(\alpha_{tx^{-1} xs^{-1}}(i(a))) \Vtil_{tx^{-1}}\Vtil_x \xi\\
    & = & \Vtil^*_{tx^{-1}} \pi(\alpha_{ts^{-1}}(i(a))) \Vtil_{tx^{-1}}\Vtil_x \xi\\
    & = & \Vtil^*_{x^{-1}}\Vtil^*_{t} \pi(\alpha_{ts^{-1}}(i(a))) \Vtil_t  \xi\\
    & = & \Vtil_{x} \pitil(\beta^{-1}_{s}(i(a))) \xi,
\end{eqnarray*}
and since $\mathcal H_0$ is dense in $\mathcal H$ and $B_0$ is dense in $B$,
the pair $(\pitil, \Vtil)$ satisfies the covariance relation.
\end{proof}

\subsection{Full corners.}
Once we know how to dilate covariant representations from the semigroup action to
the group action we can establish the relation between the respective crossed products.
Before proving our main result we recall that if $p$ is a projection in the C*-algebra $A$
then the algebra $p A p$ is  a {\em corner} in $A$, which is said to be {\em full} if 
the linear span of $ApA$ is dense in $A$.  The most relevant feature of full corners is that
if $pAp$ is a full corner in $A$, then $pA$ is a full 
Hilbert bimodule implementing the Morita equivalence, in the sense of Rieffel, of $pAp$ to $A$. 

\begin{theorem}\label{fulcor}
Suppose $(A,S,\alpha, \lambda)$ is a twisted semigroup dynamical system 
in which $S$ is an Ore semigroup acting by injective endomorphisms 
and  $\lambda $ is a multiplier on $S$.
Let $(B,G,\beta, \mu)$ be the minimal automorphic dilation,
with embedding $i: A \to B$.
Then  $A\rtimes_{\alpha,\lambda} S$ is canonically isomorphic to
 $i(1) (B\rtimes_{\beta, \mu} G) i(1)$, which is a full corner.
As a consequence, the crossed product $A\rtimes_{\alpha,\lambda} S$ is 
Morita equivalent to $B \rtimes_{\beta, \mu} G$.
\end{theorem}

\begin{proof}
Let $U$ be the projective unitary representation of $G$ in the multiplier 
algebra of $B\rtimes_{\beta, \mu} G$, and notice that 
$$i(1) U_s i(1) = U_s i(1), \qquad s \in S, $$ because $i(A)$ is invariant under 
$\beta_s$. Define $v_s =  U_s i(1) $.
Then $v_s^* v_s = i(1) U_s ^* U_s i(1) = i(1)$  and $v_s v_t = U_s i(1) U_t i(1) = 
U_s U_t i(1) = \mu(s,t) U_{st} i(1) = \lambda(s,t) v_{st}$, so $v$ is a projective isometric 
representation of $S$ with multiplier $\lambda$. Since $i(1) (B\rtimes_{\beta, \mu} G) i(1)$
is generated by the elements $i(1)U_x^* i(a) U_y i(1) = v_x^* i(a) v_y$, the
isomorphism will be established by uniqueness of the crossed product once we  
 show that the pair $(i,v)$ is universal.

Suppose $(\pi,V)$ is a covariant representation for the twisted system 
$(A,S,\alpha,\lambda)$, and let $(\pitil,\Vtil)$ be the  
corresponding dilated covariant representation of $(B,G,\beta, \mu)$
given by \lemref{dil-cov-rep}.
By the universal property of $ B\rtimes_{\beta, \mu} G $
there is a homomorphism 
$$
(\pitil \times \Vtil) : B\rtimes_{\beta, \mu} G \to C^*(\pitil, \Vtil)
$$
such that 
$\pitil(b) \Vtil_s = (\pitil \times \Vtil) (i_B(b) U_s)$ .

Let $\rho $ be the restriction of $ (\pitil \times \Vtil) $ to 
$i(1) (B\rtimes_{\beta, \mu} G) i(1)$, cut down to the invariant subspace $H$.
 By \lemref{dil-cov-rep}
$$
\rho(i(a)) = (\pitil \times \Vtil)(i(a)) = \pitil \circ i(a) = \pi(a), \qquad a \in A, 
$$
while
$$
\rho(v_s) = (\pitil \times \Vtil) (U_s i(1))  = \Vtil_s \pi(1) = V_s, \qquad s \in S
$$
Thus $\rho \circ i = \pi$ and $\rho \circ v = V$, so $(i,v)$ is universal for
$(A,S,\alpha, \lambda)$.

Finally  we prove that the corner is full, i.e., that the linear span of the elements of the form 
$X i(1) Y$ with $X, Y \in B \rtimes_{\beta, \mu} G $ is a dense subset of
 $ B \rtimes_{\beta, \mu} G$. 
It is easy to see
that the elements of the form $U_s^* b U_t$ span a dense subset of 
$ B \rtimes_{\beta, \mu} G$ because $G = S^{-1} S$, where $b$ may be replaced with $U_r^* i(a) U_r$ 
by minimality of the dilation.  
Thus  the elements $U_y^* i(\alpha_z(a)) U_x$ with $x,y,z \in S$
and $a \in A$ span a dense subset of $B \rtimes_{\beta, \mu} G$,  and since
 $i(\alpha_z(a)) = i(1) i(\alpha_z(a))$,  the proof is finished.
\end{proof}

\begin{remark} \label{nontriv}
If one drops the assumption of injectivity of the endomorphisms, it is still possible to
carry out the constructions and the arguments in the proofs of the preceding theorems. 
However, the resulting homomorphism $i: A \to B$ may not 
be an embedding any more. Indeed,  Example 2.1(a) of \cite{sta} shows that the limit algebra
$B$ may turn out to be the $0$ C*-algebra, yielding a trivial dilated system.

We notice  that  the dilated system $(B,G,\beta,\mu)$  has
nontrivial covariant representations if and only if $B \neq 0$, and these 
representations, when 
cut down to $i(A)$, give nontrivial covariant representations of the
original semigroup system $(A,S,\alpha,\lambda)$.  Thus, following
\cite[Proposition 2.2]{sta} which deals with the case $S = \N$,
we conclude  
that the crossed product $A\rtimes_{\alpha,\lambda} S$ is nontrivial if and only if 
the limit algebra $B$ is not $0$. 
Clearly, this is the case when, for instance, the endomorphisms are injective.
\end{remark}

\section{An example from number theory.}

As an application of the preceding theory we consider 
the semigroup dynamical system from \cite{bc-alg} whose crossed product 
is the Bost-Connes Hecke C*-algebra \cite{bos-con}. 
Since Morita equivalence implies that the representation theory of the semigroup 
dynamical system is equivalent to that of the dilated system, it is
quite useful to have an explicit formulation of the dilation.
We point out that since the semigroup in question is abelian,
this application is somewhat independent from the rest of the material
on nonabelian semigroups. 
In fact, the example could be dealt with by  enhancing
\cite[Section 2]{murnew} with the uniqueness and fullness properties
discussed above, which are easier to prove for abelian semigroups.

\subsection{Finite Adeles.}
The natural setting for identifying the ingredients of the minimal automorphic dilation
of the semigroup dynamical system introduced in \cite{bc-alg} 
will be the (dual)  $p$-adic picture described in \cite[Proposition 32]{diri}, in which 
the algebra is $C(\prod_p \Zp)$ and the endomorphisms $\alpha_n$ 
consist of `division by $n$' in $\prod_p \Zp$: 
$$\alpha_n (f) (x) = \left\{\begin{array}{ll}
f(x/n)&\mbox{if $n | x$}\\
0&\mbox{ otherwise}.
\end{array} \right.
$$
By \cite[Corollary 2.10]{bc-alg} the crossed product associated to this system 
is canonically isomorphic to the Bost-Connes Hecke C*-algebra $\mathcal C_{\Q}$.

The ring ${\mathcal Z} := \prod_p \Zp$ has lots of zero divisors and hence no fraction field. 
However, the diagonally embedded copy of $\nat$ 
is a multiplicative set with no zero divisors, and
we may enlarge ${\mathcal Z}$ to a ring in which division by an element of $\nat$ is
always possible. Our motivation is to extend the endomorphisms $\alpha_n$ 
defined above to automorphisms.

The algebraic part is easy: we consider the ring $(\nat)^{-1} {\mathcal Z}$ of formal fractions
${z}/{n}$ with $z\in {\mathcal Z}$ and $n \in \nat$,
 with the obvious rules of addition and multiplication (and simplification!), 
\cite[II.\S3]{lan-alg}. This ring has a universal property with respect to homomorphisms 
of ${\mathcal Z}$ that send $\nat$ into units.
Since $\nat $ has no zero divisors, the canonical
map $z  \rightarrow {z}/{1}$ is an embedding of ${\mathcal Z}$ in $(\nat)^{-1} {\mathcal Z}$.

The topological aspect requires a moments thought, after which we declare that the subring
${\mathcal Z}$ must retain its compact topology and be relatively open.
Since we want division by $n\in \nat$ to be an automorphism, this
determines a topology on the compact open sets 
$(1/n) {\mathcal Z}$ and hence on their union,
 $(\nat)^{-1} {\mathcal Z}$,
which becomes a locally compact ring containing ${\mathcal Z}$ as a compact open subring.

The ring we have just defined is (isomorphic to) the locally compact ring $\mathbb A_f$ of
finite adeles, which is usually defined as the restricted product, over the primes $p \in \primes$
of the $p$-adic numbers $\Q_p$ with respect to the $p$-adic integers $\Zp$:
$$
\mathbb A_f : = \{  (a_p) : a_p \in\mathbb Q_p \text{ and } a_p \in \Zp \text{ for all but finitely
many }p \in \primes \},
$$
with $ \prod_p \Zp$ as its maximal compact open subring.
The isomorphism is implemented by the map from $(\nat)^{-1} {\mathcal Z}$ into $\mathbb A_f$
given by the universal property; this map is clearly injective and, 
since every finite adele can be written as $z/n$ with $z\in {\mathcal Z}$ and $n \in \nat$,
it is also surjective.
Specifically, for each $a_p \in \Q_p$ there exists $k_p$ such that $p^{k_p} a_p = z_p \in \Zp$
 and a sequence $a = (a_p)_{p\in \mathcal P}$ is an adele 
if and only if $k_p$ can be taken to be $0$ for all but finitely many
$p$'s, in which case $ n = \prod_p p^{-k_p} \in \nat$ and $a = (na)/n$,
 with $na = (na_p)_{p\in \primes} \in \prod_p \Zp$.

\subsection{The minimal automorphic dilation of $(C(\mathcal Z), \nat, \alpha)$.}\label{diladeles}
The rational numbers are embedded in $\mathbb A_f$, and division by a nonzero rational is 
clearly a homeomorphism so
$$
\beta_r(f)(a) = f(r^{-1} a), \quad a \in\mathbb A_f, r \in \qat 
$$
defines an action of $\qat = (\nat)^{-1} \nat $ by automorphisms of $C_0(\mathbb A_f)$.
 
Since ${\mathcal Z}$ is compact and open, 
its characteristic function $ 1_{{\mathcal Z}}$ is a projection in 
$C_0(\mathbb A_f)$ and there is an obvious embedding $i$ of $C({\mathcal Z})$
as the corresponding ideal of $ C_0(\mathbb A_f)$,  given by 
$$
i(f) (a) = \left\{ \begin{array}{ll}
                     f(a) & \text{ if }  a \in {\mathcal Z}\\
                       0   &  \text{ if } a \notin {\mathcal Z}.
                         \end{array} \right.
$$

\begin{proposition}\label{BC-min-aut-ext}
The C*-dynamical system $(C_0(\mathbb A_f), \qat, \beta)$ is the minimal
automorphic dilation of the semigroup dynamical system 
$(C({\mathcal Z}), \nat, \alpha)$, and hence $\mathcal C_{\Q}$ is the full corner of 
$C_0(\mathbb A_f)\rtimes_\beta \qat$ determined by the projection $1_{\mathcal Z}$. 
\end{proposition}

\begin{proof}
The embedding clearly intertwines $\alpha_n$ and $\beta_n$, in the sense that
$\beta_n (i(f))  = i(\alpha_n(f))$, and 
the union of the compact subgroups $(1/n) {\mathcal Z}$
is dense in $\mathbb A_f$, so the union of the
subalgebras $\beta_{1/n}(i( C({\mathcal Z})) )$
is dense in $C_0(\mathbb A_f) $, and the result follows from 
\thmref{dil-ext} and \thmref{fulcor}. 
\end{proof}

\medskip
Since the discrete multiplicative group $\qat$ acts by homotheties on the
locally compact additive group $\mathbb A_f$, and since $\mathbb A_f$ is self-dual,
we obtain another characterization of $\mathcal C_{\Q}$ as a full corner in the 
group C*-algebra of the semidirect product $\mathbb A_f \rtimes \qat$.
One should bear in mind, however, that the self duality of $\mathbb A_f$ is not canonical.

\begin{corollary}
Let $e_{{\mathcal Z}} \in C^*(\mathbb A_f)$ 
be the Fourier transform of ${1_{\mathcal Z}} \in C_0(\mathbb A_f)$. Then  
$$ \mathcal C_{\Q} \cong e_{{\mathcal Z}} C^*(\mathbb A_f \rtimes \qat)e_{{\mathcal Z}}.$$
\end{corollary}
\begin{proof}
The action of $\qat$ on $\mathbb A_f$ is by 
homotheties, which are group automorphisms, so
$C^*(\mathbb A_f \rtimes \qat)$ is isomorphic to
the crossed product $C^*(\mathbb A_f) \rtimes_\beta \qat$.
Moreover, the self-duality of the additive group of $\mathbb A_f$ satisfies
$\langle rx,y\rangle = \langle  x,ry\rangle$ for $r \in \qat$,   
thus  $C^*(\mathbb A_f)$ is covariantly isomorphic to $C_0(\mathbb A_f)$, so
 $C^*(\mathbb A_f) \rtimes_\beta \qat$ is isomorphic 
to $C_0(\mathbb A_f) \rtimes_\beta \qat$, 
and the claim follows from Proposition \ref{BC-min-aut-ext}. 
\end{proof} 

\begin{remark}
One of the principles of noncommutative geometry
advocates that if $G$ is a group acting on a space $X$, then the quotient space $X/G$ 
has a noncommutative version 
in the associated crossed product $C_0(X) \rtimes G$, which is often more tractable.
Accordingly, if we allow back in the all-important
place at infinity which is left out from $\mathcal A_f$ 
and if we substitute $\qat$ by $\Q^*$, cf. \cite[Remarks 33]{bos-con},
then our \proref{BC-min-aut-ext} gives an explicit path leading from
the  Bost-Connes Hecke  C*-algebra to the  space $\mathcal A/\mathbb Q^*$,
on which the construction of \cite{con-cr, con-rzf} is based. 
\end{remark}

\end{document}